\newcommand{\myaddress}{paul@math.tau.ac.il}

\newcommand{\SIZE}{11pt}

\newcommand{\LTX}{\documentstyle{article}
                  \newcommand{\Bbb}{\bf}
                  \newcommand{\myaddress}{paul@math.tau.ac.il}
                 }

\newcommand{\AMSLTX}{\documentstyle[\SIZE , amscd, amssymb]{amsart}}

\AMSLTX

\newtheorem{thm}{Theorem}[section]
\newtheorem{lem}[thm]{Lemma}

\newtheorem{prop}[thm]{Proposition}

\newtheorem{cor}[thm]{Corollary}

\newtheorem{cnj}[thm]{Conjecture}

\newcommand{\Qed}{\hfill \rule{.75em}{.75em}}

\newcommand{\cntrs}{\setcounter{thm}{0} \renewcommand{\thethm}
{\thesection.\Alph{thm}}}
\newcommand{\cntrsb}{\setcounter{thm}{0} \renewcommand{\thethm}
{\thesubsection.\Alph{thm}}}

\newcommand{\PoD}{Poincar\'{e} dual }
\newcommand{\Khlr}{K\"{a}hler }
\newcommand{\rarw}{\rightarrow}
\newcommand{\lrarw}{\longrightarrow}

\newcommand{\stmin}{\setminus}
\newcommand{\lbr}{\langle}   		
\newcommand{\rbr}{\rangle}   		
\newcommand{\qd}{\quad}

\newcommand {\CPTU}{{{\Bbb C}P^2}}

\newcommand {\om}{\omega}
\newcommand {\omstd}{\omega_{std}}
\newcommand {\Om}{\Omega}

\newcommand {\sig}{\sigma}
\newcommand {\Sig}{\Sigma}

\newcommand {\lam}{\lambda}

\newcommand {\eps}{\epsilon}
\newcommand {\dlt}{\delta}
\newcommand {\al}{\alpha}

\newcommand {\Th}{\Theta}
\newcommand {\vphi}{\varphi}

\newcommand {\CC}{{\Bbb C}}
\newcommand {\NN}{{\Bbb N}}
\newcommand {\RR}{{\Bbb R}}
\newcommand {\QQ}{{\Bbb Q}}
\newcommand {\ZZ}{{\Bbb Z}}

\newcommand {\PP}{{\Bbb P}}

\newcommand {\calL}{{\cal{L}}}
\newcommand {\calC}{{\cal{C}}}
\newcommand {\calB}{{\cal{B}}}

\newcommand {\calF}{{\cal{F}}}
\newcommand {\calO}{{\cal{O}}}
\newcommand {\calE}{{\cal{E}}}
\newcommand {\calR}{{\cal{R}}}
\newcommand {\calK}{{\cal{K}}}
\newcommand {\calJ}{{\cal{J}}}
\newcommand {\calP}{{\cal{P}}}
\newcommand {\calU}{{\cal{U}}}
\newcommand {\calV}{{\cal{V}}}

\newcommand {\calKbar}{\overline{\calK}}
\newcommand {\wtld}{\widetilde}
\newcommand {\wtldOm}{\widetilde{\Om}}

\newcommand {\wtldM}{{\widetilde{M}}}

\newcommand {\wtldS}{{\widetilde{S}}}
\newcommand {\wtldD}{{\widetilde{D}}}
\newcommand {\wtldC}{{\widetilde{C}}}
\newcommand {\wtldX}{{\widetilde{X}}}

\newcommand {\undl}{\underline}

\newcommand {\phipack}{\vphi=\coprod_{j=1}^N
\vphi_j : \coprod_{j=1}^N B(\lam_j) \rarw (M,\Om)}

\newcommand {\Oone}{{\cal O}_{\CPTU}(1)}
\newcommand {\dzi}{{\#_{_i}}}

\begin{document}

\title{Constructing new ample divisors out of old ones} 
\date{May 1997}
\author{Paul Biran}
\address{School of Mathematics\\ 
Tel-Aviv University\\
Ramat-Aviv \  Tel-Aviv 69978\\ 
Israel.}
\email{\myaddress}
\thanks{\bf Preliminary Version.}

\begin{abstract}

We prove a gluing theorem which allows to construct an ample divisor on a 
rational surface from two given ample divisors on simpler surfaces. 
This theorem combined with the Cremona action on the ample cone gives rise to 
an algorithm for constructing new ample divisors. 
We then propose a conjecture relating continued fractions 
approximations and Seshadri-like constants of line bundles over rational 
surfaces. By applying our algorithm recursively we verify our conjecture in 
many cases and obtain new asymptotic estimates on these constants.
Finally, we explain the intuition behind the gluing theorem in terms of 
symplectic geometry and propose generalizations. 

\end{abstract}

\maketitle

\section{Introduction} \cntrs
\label{sect-Intro}

The main objective of this paper is to propose a method for constructing new
ample divisors on rational surfaces by gluing two given ones.

Recall that a divisor $D$ on an algebraic variety $X$ is ample if the 
corresponding line bundle
$\calO_X(D)$ is ample, and is called {\em nef (numerically effective)}
if there exists an ample divisor $A$ such that $A+kD$ is ample for every
$k>0$. We refer the reader to \cite{Dem,Ha-Ample} for excellent expositions 
on various aspects of the theory of ample and nef line bundles.

Of fundamental importance is the determination of those classes in 
$\mbox{Pic}(X)$ 
which are ample. Although this problem has a very simple solution for smooth 
curves, already in dimension two the problem becomes much harder. 
It turns out the even for relatively simple surfaces,
such as rational, the complete answer is not known. Several conjectures in
this direction exist, however at the present time only estimates on the
{\em ample cone} -- the cone generated by the ample classes in 
$\mbox{Pic}(X)$ -- are known. 
For example, let $d,m>0$ and consider the divisor class
$$D=\pi^* \calO_{\CPTU}(d) - m\sum_{j=1}^N E_j$$ on the blow-up 
$\pi:V_N \rarw \CPTU$ of $\CPTU$ at $N\geq 9$ generic points. 
Nagata conjectured in~\cite{Nag} that {\em $D$ is ample iff $D\cdot D >0$}, 
but was able to prove it only for $N's$ which are squares. 
In~\cite{Xu-Curves} Xu proved that {\em $D$ is ample provided that 
$\frac{m}{d}< \frac{\sqrt{N-1}}{N}$}.
By making a more detailed analysis of the case $m=1$, Xu proved in
\cite{Xu-Divisors} that {\em when $d\geq 3$ the divisor class 
$D=\pi^* \calO_{\CPTU}(d) - \sum_{j=1}^N E_j$ is ample iff $D \cdot D > 0$} 
(see also K\"{u}chle~\cite{Ku} for a generalization for arbitrary surfaces 
and~\cite{Ang} for an analogous result for $\CC P^3$).

Closely related is the problem of computing {\em Seshadri constants}
of ample line bundles, which measure their local positivity.
The Seshadri constant $\calE(\calL,p)$ of the line bundle $\calL$ at
the point $p\in X$ is defined to be {\em the supremum of all those 
$\eps\geq 0$ for which the $\RR$-divisor class $\pi^*\calL -\eps E$ is nef on 
the blow-up $\pi:\wtldX_p \rarw X$ of $X$ at the point $p$ with exceptional 
divisor $E$}.

Seshadri constant has been studied much by Demailly (\cite{Dem}),
Ein, K\"{u}chle, Lazarsfeld (\cite{EL},\cite{EKL},\cite{Laz}), 
and Xu (\cite{Xu-Ample}). 
A considerable part of these works is devoted to computations
and estimates from below on the values of these constants. 

The present paper is largely motivated by the problem of computing Seshadri
constants and the determination of the ample cone of rational surfaces.
Our main results provide an algorithmic method for constructing new ample
divisors out of the knowledge of ample divisors on simpler rational surfaces.
By applying the algorithm recursively we obtain in 
Section~\ref{sect-Asymptotics} new estimates on Seshadri-like constants and 
detect new ample divisors.
We then propose in Section~\ref{sect-Conjecture} a conjecture naturally 
arising from our method which relates {\em continued fractions expansions of 
$\sqrt{N}$} with the ample cone of $\CPTU$ blown-up at $N$ points.
Finally we interpret in Section~\ref{sect-Symplectic} our main results in 
the language of Symplectic Geometry and explain the intuition behind them.

Our main tool is Shustin's version of the Viro method for gluing curves
with singularities.

\section{Main results} \cntrs

Our main results deal with {\em simple rational surfaces} $S$, which by
definition are blow-ups $\Th:S\rarw \CPTU$ of $\CPTU$ at $n$ 
distinct points $p_1, \ldots, p_n \in \CPTU$.\footnote{Note that we regard 
$\CPTU$ itself as a simple rational surface too (this corresponds to $n=0$).} 
We denote by $E^S_i=\Th^{-1}(p_i) \qd i=1,\ldots,n$ the 
{\em standard exceptional divisors} of the blow-up and write $\Sig^S$ for 
the union $\cup_{i=1}^n E^S_i$. 
Finally, we write $L^S$ for be a divisor on $S$, obtained by pulling back 
via $\Th$ a projective line in $\CPTU$ which 
does not pass through any of the points $p_1,\ldots,p_n$. 

A vector 
$(d;\al_1,\ldots,\al_k) \in \ZZ_+ \times \ZZ^k_{\mbox{\tiny $\geq 0$}}$ is
called ample (resp. nef) if there exists a simple rational surface $V$,
on which the divisor $dL^V-\sum_{j=1}^k \al_j E^V_j$ is ample
(resp. nef). \ \\ 

Our first result is the following gluing theorem:
\begin{thm}
\label{thm-glue1}

Let $(d;m_1, \ldots m_n, m)$ be an ample (resp. nef) vector and 
$(m; \al_1, \ldots, \al_k) \in \ZZ_+^{k+1}$ a nef vector. Then
$v=(d;m_1,\ldots,m_n,\al_1,\ldots, \al_k)$ is ample (resp. nef). 
Moreover, $v$ can be realized by an ample (resp. nef) divisor on a very 
general rational surface.

\end{thm}
By a {\em very general} choice of points $q_1,\ldots,q_r$ in an
algebraic variety $X$ we mean that $(q_1,\ldots ,q_r)$ is allowed to
vary in a subset of the configuration space
$\calC_r(X)=\{(x_1, \ldots, x_r)\in X^r \mid x_i \neq x_j \}$
whose complement is contained in a countable union of proper
subvarieties of $\calC_r(X)$. By a {\em very general rational surface}
we mean one which is obtained by blowing-up points $q_1,\ldots,q_r \in
\CPTU$ which may be chosen to be very general.

We shall actually prove a stronger result which allows us to keep the
blown-up points corresponding to the first ample vector fixed, thus giving 
information also on ample divisors on non-generic rational surfaces. 
The precise statement is:

\begin{thm}
\label{thm-glue2}

Let $D$ be a divisor on a simple rational surface $S$. Suppose that
there exists a point $p\in S\setminus (\Sig^S \cup Supp\,D)$ and
$m>0$ such that $\pi_p^* D-mE$ is ample on the blow-up
$\pi_p:\wtldS_p\rarw S$ of $S$ at $p$ with exceptional divisor $E$.
Let $(m;\al_1, \ldots, \al_k) \in \ZZ_+^{k+1}$ be a nef vector. Then for
a very general choice of points $q_1, \ldots, q_k \in
S\setminus(\Sig^S\cup Supp\,D)$ the divisor
$$\pi^*D -\sum_{j=1}^k \al_j E_j$$ is ample on the blow-up 
$\pi:\wtldS\rarw S$ of $S$ at $q_1, \ldots, q_k$ with exceptional divisors 
$E_j=\pi^{-1}(q_j)$.

\end{thm}
Proofs of Theorems~\ref{thm-glue1} and~\ref{thm-glue2} appear in 
Section~\ref{sect-Gluing}. 

Theorem~\ref{thm-glue1} in combination with the action of the Cremona
group on the ample cone give rise to an algorithmic procedure for detecting 
new ample classes in the Picard group of rational surfaces. The algorithm
will be explained in Section~\ref{subsect-Alg}. 

\subsection{Applications to Seshadri constants} \cntrsb
\label{subsect-Applic}

Given an ample line bundle $\calL \rarw S$ on a surface, and a vector
$w=(w_1,\ldots, w_N)$ of positive numbers we define the {\em
$w$-weighted remainder} of $\calL$ at the $N$ distinct points
$p_1,\ldots,p_N \in S$ to be the quantity
$$\calR^w(\calL, p_1,\ldots,p_N)=\frac{1}{\calL \cdot \calL}
\inf_{0\leq\eps\in \RR} \left\{ \calL_{\eps} \cdot \calL_{\eps} \bigg| 
\calL_{\eps} = \pi^*\calL -\eps\sum_{j=1}^N w_j E_j \quad \mbox{is nef} 
\right\},$$
where $\pi:\wtldS\rarw S$ is the blow-up of $S$ at the points
$p_1,\ldots,p_N$ with exceptional divisors $E_i=\pi^{-1}(p_i)$.
It is obvious that $0 \leq \calR^w < 1$.
Note that $\calR^w$ remains invariant
under rescalings of $\calL$ and of $w$, namely $\calR^{aw}(b\calL, p_1,
\ldots p_N)=\calR^{w}(\calL, p_1,\ldots p_N)$ for every $a,b>0$.
It is convenient to define also a more global invariant, namely
$$\calR^w_N (\calL)= \inf \left\{ \calR^{w}(\calL, p_1,\ldots p_N) 
\mid p_1, \ldots p_N \in S \quad \mbox{are distinct points} \right\}.$$

Restricting to the case of homogeneous weights we obtain the 
{\em homogeneous remainders}
$$\calR(\calL, p_1, \ldots, p_N)= \calR^{w_h}(\calL, p_1, \ldots, p_N), 
\qquad  \calR_N(\calL) = \calR^{w_h}(\calL),$$
where $w_h=(1, \ldots, 1)$.

The constants $\calR^w(\calL,p_1, \ldots,p_N)$ are obvious generalizations of 
the Seshadri constants $\calE(\calL,p)$ from section~\ref{sect-Intro}
(see also~\cite{Xu-Ample} for similar Seshadri-like constants). 
Several theorems and conjectures related to the ample cone can be neatly
formulated using the constants $\calR_N$. For example, Nagata's
conjecture from Section~\ref{sect-Intro} can be reformulated as
``{\em \ $\calR_N(\Oone)=0$ when $N \geq 9$ \ }''.
Similarly, Xu's result from Section~\ref{sect-Intro} asserts that 
$\calR_N(\Oone) \leq \frac{1}{N}$. In Section~\ref{subsect-prfs1} we shall 
prove the following asymptotic result:

\begin{thm}
\label{thm-asymp1} 
\begin{enumerate}
\item[1)] For $N=a^2l^2+2l$, $a,l\in \NN$ \ \ $\calR_N(\Oone)\leq
\frac{1}{(a^2l+1)^2}$.  
\item[2)] For $N=a^2l^2-2l$, $a,l,\in \NN$ \ \  $\calR_N(\Oone)\leq
\frac{1}{(a^2l-1)^2}$. 
\item[3)] If $N=a^2l^2+l$ with $1<l\in \NN$ and suppose that 
$l > \frac{a}{2^{k-1}}$, where $k$ is the maximal non-negative integer for 
which $a\equiv 0 \mod 2^k$. Then $\calR_N(\Oone)\leq\frac{1}{(2a^2l+1)^2}$.
\end{enumerate}
\end{thm}
In Section~\ref{sect-Conjecture} we shall view this result in a more general 
context by proposing a conjecture which bounds $\calR_N(\Oone)$ in terms of 
continued fractions approximations of $\sqrt{N}$.

Our methods also yield, as a corollary, the following generalization of a 
theorem of Xu~\cite{Xu-Divisors} and K\"{u}chle~\cite{Ku}:
\begin{cor}
\label{cor-coef2}

Let $d>0$. The divisor $D=\pi^* \calO_{\CPTU}(d)-2\sum_{j=1}^N E_j$ on the
blow-up of $\CPTU$ at $N$ very general points is nef iff $D\cdot D\geq 0$.


\end{cor}
The proof appears in Section~\ref{subsect-prfs1}.
Let us conclude this section with the following, somewhat amusing,
corollary of Theorem~\ref{thm-glue1}.
\begin{cor}
\label{cor-Nag}
If Nagata's conjecture holds for $N_1$ and $N_2$ then it holds also for 
$N_1N_2$.
\end{cor}
The proof is given in Section~\ref{subsect-prfs2}.

\section{Gluing curves on rational surfaces} \cntrs
\label{sect-Gluing}

We shall derive Theorem~\ref{thm-glue1} as a corollary from
Theorem~\ref{thm-glue2}. The proof of Theorem~\ref{thm-glue2} is
based on a technique for "gluing" singular curves, which was developed
by Shustin in~\cite{Sh1}. This method generalizes Viro's method 
(see~\cite{Vir}) for gluing curves to singular cases. 

Suppose that $C_1, \ldots, C_n$ are plane curves with
Newton Polygons $\Delta_1, \ldots, \Delta_n$ which have mutually
disjoint interiors and match together to a bigger polygon $\Delta =
\Delta_1 \cup \ldots \cup \Delta_n$. Shustin's method allows,
under some transversality conditions on the equisingular strarta
corresponding to $C_1, \ldots, C_n$, to construct a new curve $C$ with
Newton polygon $\Delta$ and with singular points "inherited" from
$C_1,\ldots,C_n$. We refer the reader to~\cite{Sh1} for a detailed 
presentation of the general method and to~\cite{Sh2} for interesting 
applications in other directions. Here, we shall make use only of a tip of 
the power of this method, by applying it to two curves with disjoint Newton 
polygons. 

The application of Shustin's technique to our problem is summed up in the
following proposition which will be the main ingredient in the proof of
Theorem~\ref{thm-glue2}. Most of the proof presented below is essentially an 
adjustment of the arguments appearing in the proof of Theorem 3.1 
of~\cite{Sh1} to our specific situation. 

\begin{prop}
\label{prop-glue}

Let $D$ be an effective divisor on a simple rational surface $S$, and
$p\in S \setminus\Sig^S$ a point with $\mbox{\em mult}_p D=m>0$.
Let $C$ be an effective divisor on another simple rational surface $V$,
lying in the linear system $|m' L^V-\sum_{j=1}^k \al_j E^V_j|$.
Suppose that $D,C$ satisfy the following conditions: 
\begin{enumerate}
\item[1)] $0<m'<m$. 
\item[2)] $H^1(\wtldS_p, \calO_{\wtldS_p}(\pi^*_p D-mE))=0$, where
$\pi_p:\wtldS_p\rarw S$ is the blow-up of $S$ at the point $p$ with
exceptional divisor $E$. 
\item[3)] $H^1(V,\calO_V(C))=0$. 
\item[4)] Each of $C,D$ does not have any of the standard exceptional divisors
$E^V_i, E^S_j$ as one of its components. 
\item[5)] $D$ is an irreducible curve. 
\end{enumerate}
Then, there exist $k$ distinct points
$q_1, \ldots q_k \in S \setminus (\Sig^S \cup Supp\, D)$ and a curve
$\wtld{D}$ on the blow-up $\pi:\wtldS \rarw S$ of $S$ at
$q_1,\ldots,q_k$ with exceptional divisors $E_j=\pi^{-1}(q_j)$, which has the 
following properties: 
\begin{enumerate}
\item[1)] $\wtld{D}\in |\pi^*D-\sum_{j=1}^k \al_j E_j|$. 
\item[2)] The curve $\wtld{D}$ is irreducible. 
\end{enumerate}
\end{prop}

\pf The idea of the proof is basically the following. By passing to the 
underlying projective planes of $\wtldS_p$ and $V$ we obtain from $D$ and $C$ 
two singular curves $C_1$ and $C_2$ and a point, still denoted by $p$, 
such that $\mbox{mult}_p C_1 > \deg C_2$. This inequality implies that 
the Newton polygons of $C_1$ and $C_2$ with respect to an affine chart 
centered at $p$ are disjoint. 
The next step is to construct two deformations $C_{1,t}$ 
and $C_{2,t}$ of $C_1$ and $C_2$ which are equisingular for $t>0$ and such 
that each of them contains a deformations of the union of the singular  
points of $C_1$ and $C_2$ except of the one at the point $p$ which might 
disappear. These two deformations are then glued using the 
{\em Viro polynomial}. Shustin's method requires the Newton polygons of each 
of $C_{1,t}$ and $C_{2,t}$ to be contained in the union, say $\Delta$, of the 
ones of $C_1$ and $C_2$. In order to construct deformations which satisfy 
this, one has to prove roughly speaking that the equisingular strata of 
$C_1$ and $C_2$ intersect transversally the space of curves with 
Newton polygons $\Delta$. This is precisely what the conditions of vanishing 
of the $H^1$'s is needed for. 
 
Let us give now the precise details of the proof.
Suppose that $S$ is obtained by blowing-up $\Th_{S}:S \rarw \CPTU$ at 
$p_1, \ldots, p_N \in \CPTU$ and that $V$ is obtained by blowing-up 
$\Th_{V}:V \rarw \CPTU$ at $q_1^0, \ldots, q_k^0 \in \CPTU$.
Put $p_0=\Th_{S}(p), \quad C_1=\Th_{S}(D)\subset \CPTU$ and 
$C_2=\Th_{V}(C)\subset \CPTU$. 
Assuming that $D\in|dL^S-\sum_{i=1}^n m_i E^S_i|$ we see that:

\begin{itemize}
\item  $C_1$ is a plane curve of degree $d$ and has singularities of orders 
$m_1, \ldots, m_n$ at the points $p_1, \ldots, p_n$ and a singular point of 
order $m$ at $p_0$. 
\item $C_2$ is a plane curve of degree $m'$ and has singularities of orders 
$\al_1,\ldots,\al_k$ at the points $q_1^0, \ldots, q_k^0$.
\end{itemize}
In view of what we have to prove there is no loss of generality in assuming 
that $q_j^0 \not\in C_1$ for every $1\leq j\leq k$. 
Choose an affine chart $\CC^2 \subset \CPTU$ with coordinates $(x,y)$ such 
that $p_0=(0,0)$ and such that 
$p_1,\ldots, p_n, q_1^0, \ldots, q_k^0 \in (\CC^*)^2 \subset \CC^2 
\subset \CPTU$. 

Let $F_1(x,y),F_2(x,y)$ be polynomials of degrees $d$ and $m'$ 
respectively, such that $C_1 \cap \CC^2 = \{ F_1=0\}$ and 
$C_2\cap \CC^2 = \{ F_2=0 \}$.
Set $$\Delta_1=\{(i,j) \in \ZZ^2_{\geq 0} | m\leq i+j \leq d\}, \qquad 
\Delta_2 = \{(i,j) \in \ZZ^2_{\geq 0} | 0\leq i+j \leq m' \},$$
and put $\Delta = \Delta_1 \cup \Delta_2$. With these notations, we may write 
$$F_1(x,y)=\sum_{(i,j)\in \Delta_1} a_{ij}x^i y^j, \qquad
  F_2(x,y)=\sum_{(i,j)\in \Delta_2} a_{ij}x^i y^j.$$
Let $\overline{\Delta}_1\supset \Delta_1$ and 
$\overline{\Delta}_2 \supset \Delta_2$ be two slightly larger triangles with 
disjoint interiors. More precisely, let $\dlt>0$ be a small enough number 
such that $m'+2\dlt < d-2\dlt$ and set 
$$ \overline{\Delta}_1=\{(i,j) \in \RR^2_{\geq 0} | m-\dlt \leq i+j \leq d\},
\qquad
\overline{\Delta}_1=\{(i,j) \in \RR^2_{\geq 0} | 0\leq i+j \leq m'+\dlt \}. $$
Next, choose a strictly convex continuous piecewise linear function 
$\nu:\RR^2 \rarw \RR$ such that the restrictions of $\nu$ to each of 
$\overline{\Delta}_1,\overline{\Delta}_2$ coincides with some linear function 
$\ell_1,\ell_2:\RR^2 \rarw \RR$ with $\ell_1 \neq \ell_2$.
We shall define the curve $\wtld{D}$ as the zero locus of a polynomial lying 
in the following family: 
$$F_t(x,y)=\sum_{(i,j) \in \Delta} A_{ij}(t)x^i y^j t^{\nu(i,j)} 
\qquad t>0,$$ with $\lim_{t\to 0}A_{ij}(t)=a_{ij}$. The polynomial $F_t$ is 
called the {\em Viro polynomial}.

More precisely, we claim that by a correct choice of of the coefficients 
$A_{ij}(t)$, and of a homogeneous change of coordinates 
$(x,y)\rarw T_t(x,y)$, the curve 
$D_t=\{F_t(T_t(x,y))=0\}$ will have the following properties for 
$t>0$ small enough: 
\begin{enumerate}
\item[1)] $\mbox{mult}_{p_i} D_t = m_i$ for every $1\leq i \leq n$. 
\item[2)] There exists $k$ points ${q_1}_t, \ldots, {q_k}_t$ depending 
smoothly on $t>0$ such that ${q_j}_t \neq p_i$ for every $i,j$ and 
$\mbox{mult}_{{q_j}_t} D_t = \al_j$ for every $1\leq j \leq k$. 
\item[3)] $D_t$ is an irreducible curve of degree $d$. 
\end{enumerate}

If we manage to prove this then the statement of the proposition will 
immediately follow. Indeed, let $t_0>0$ be small enough such that 
properties 1-3 above hold. 
Consider $\overline{D}_{t_0} \subset \CPTU$, the closure of 
$D_{t_0}$ in $\CPTU$, and let $\wtld{D}_{t_0}$ be the proper transform of 
$\overline{D}_{t_0}$ in $\wtldS$, the blow-up of $\CPTU$ at 
$p_1,\ldots,p_n, {q_1}_{t_0}, \ldots ,{q_k}_{t_0}$. 
Clearly $\wtld{D}_t \in |\pi^*D -\sum_{j=1}^k \al_j E_j|$, where 
$\pi:\wtldS \rarw S$ denotes the blow-up of $S$ at 
$q_1={q_1}_{t_0}, \ldots ,{q_k}_{t_0}$. 

Let us prove the existence of the coefficients $A_{ij}(t)$ having the claimed
properties. For this end, set $\nu_1=\nu-\ell_1, \quad \nu_2=\nu-\ell_2$. 
Note that since $\nu$ is strictly convex and $\ell_1\neq \ell_2$ by 
construction, we must have ${\nu_1}_{|_{\overline{\Delta}_2}} > 0, \quad 
{\nu_2}_{|_{\overline{\Delta}_1}}>0$.
Consider the following deformations of $F_1(x,y), F_2(x,y)$: \\
$$F_{1,t}(x,y) = \sum_{(i,j) \in \Delta} A_{ij}(t)x^i y^j t^{\nu_1(i,j)},
\eqno(1)$$ 
$$F_{2,t}(x,y) = \sum_{(i,j) \in \Delta} A_{ij}(t)x^i y^j t^{\nu_2(i,j)}.
\eqno(2)$$
An easy computation gives:\footnote{Here we use the convention that 
$t^0 \equiv 1$ and so the families $F_{1,t},F_{2,t}$ extend smoothly to 
$t\geq 0$.} 
$$F_{1,t}(x,y) = 
F_1(x,y)+ \sum_{(i,j)\in\Delta_2} A_{ij}(t)x^i y^j t^{\nu_1(i,j)} + 
	  \sum_{(i,j)\in\Delta_1} (A_{ij}(t)-a_{ij})x^i y^j, \eqno(3)$$ 
$$F_{2,t}(x,y) = 
F_2(x,y)+ \sum_{(i,j)\in\Delta_1} A_{ij}(t)x^i y^j t^{\nu_2(i,j)} + 
	  \sum_{(i,j)\in\Delta_2} (A_{ij}(t)-a_{ij})x^i y^j, \eqno(4)$$ and 
$$F_t(x,y)=t^{c^{1}_0}F_{1,t} (t^{c^1_1}x, t^{c^1_2}y) = 
	   t^{c^2_0}F_{2,t} (t^{c^2_1}x, t^{c^2_2}y), \eqno(5)$$
where $c_i^1, c_j^2$ are the coefficients of the linear functions 
$\ell_1,\ell_2$, namely  
$$\ell_1(i,j)=c^1_0+c^1_1 i + c^1_2 j, \qquad 
\ell_2(i,j)=c^2_0+c^2_1 i + c^2_2 j.$$

Since $A_{ij}(t) \ \begin{Sb} \lrarw \\ t \to 0 \end{Sb} a_{ij}$ we see 
from (3) above that $F_{1,t} \rarw F_1$ as $t\rarw 0$. 
As $F_1(x,y)$ is assumed to be irreducible and $F_{1,t}$ is a deformation of 
$F_1$ (of the same degree) we see that $F_{1,t}$ is irreducible for $t>0$ 
small enough. In view of (5) we conclude that $F_t(x,y)$ is irreducible for 
$t>0$ small too. From (5) above we also see that the curve 
$\{ F_t(x,y)=0 \}$ will have the same (topological) types of 
singularities as each of the curves $\{F_{1,t}(x,y)=0\}, \quad 
\{F_{2,t}(x,y)=0 \}$. Put 
$$T_t (x,y)=(t^{-c^1_1}x, t^{-c^1_2}y), \qquad 
  T'_t(x,y)=(t^{-c^2_1}x, t^{-c^2_2}y)$$ and write 
$$D_t=\{F_t(\,T_t(x,y)\,)=0\}, \quad {q_j}_t=T^{-1}_t\circ T'_t (q_j^0)
\mbox{\ \ for every $1\leq j \leq k$}.$$
Clearly, the maps $(x,y)\rarw T_t(x,y)$ and $(x,y)\rarw T'_t(x,y)$ 
extend to a family of biholomorphisms of $\CPTU$ depending smoothly on $t>0$.
Note that from the definition of the points ${q_j}_t$ it easily follows that 
for a generic choice of $t>0$ the points ${q_j}_t$ will be distinct 
from the $p_i$'s. In particular there exit arbitrarily small values $t>0$ for 
which the points ${q_j}_t$ will not collide with the $p_i$'s.

Putting $D_t=\{F_t(\,T_t(x,y)\,)=0\}$, the problem is reduced to proving the 
following \\
{\bf Lemma.} There exists a smooth deformation $\{ A_{ij}(t) \}_{0<t<\eps}$ 
of the coefficients $a_{ij},\quad (i,j)\in\Delta$ with the following 
properties: 
\begin{enumerate}
\item[1)] $\lim_{t\to 0}A_{ij} = a_{ij}$. 
\item[2)] The curve $\{F_{1,t}(x,y)=0\}$ passes through 
$p_1, \ldots, p_n$ with multiplicities $m_1, \ldots, m_n$. 
\item[3)] The curve $\{F_{1,t}(x,y)=0 \}$ passes through 
$q_1^0, \ldots, q_k^0$ with multiplicities $\al_1, \ldots ,\al_k$.
\end{enumerate}

\noindent {\em Proof of the Lemma}. 
Let $\calP(\Delta_1), \calP(\Delta_2)$ be the 
spaces of polynomials in the variables $(x,y)$ with Newton diagrams contained 
in $\Delta_1,\Delta_2$, respectively. 
For every point $q\in \CC^2$, we denote by $J^{(r)}_q$ the 
space of $r$ jets of holomorphic functions at the point $q$, viewed as a 
vector space and write $j^{(r)}_q(F) \in J^{(r)}_q$ the $r$'th jet of 
$F$ at the point $q$. 

Consider the linear maps 
$$R_1:\calP(\Delta_1) \rarw \bigoplus_{i=1}^n J^{(m_i)}_{p_i}, \qquad 
  R_2:\calP(\Delta_2) \rarw \bigoplus_{j=1}^k J^{(\al_j)}_{q_j^0}$$ 
defined by  
$$R_1(F) = \left(j^{(m_1)}_{p_1}(F), \ldots, j^{(m_n)}_{p_n}(F)\right), 
\qquad R_2(F) = 
\left(j^{(\al_1)}_{q_1^0}(F), \ldots, j^{(\al_k)}_{q_k^0}(F)\right).$$ 
We claim that they are both surjective.
 
To see this let us denote for every $q\in \CPTU$ by $\calJ_q$ the ideal 
sheaf corresponding to the point $q$. 
Consider the ideal sheaf $\calJ_{X_1}=\prod_{i=1}^n \calJ^{m_i}_{p_i} 
\cdot \calJ^m_p$ on $\CPTU$, and let $X_1\subset \CPTU$ be the 
zero-dimensional subscheme defined by $\calJ_{X_1}$, with structure sheaf 
$\calO_{X_1}=\calO_{\CPTU}/\calJ_{X_1}$. 
Tensoring the structural exact sequence of $X_1$ by $\calO_{\CPTU}(d)$ 
we obtain the following exact sequence
$$0\rarw \calJ_{X_1}(d) \rarw \calO_{\CPTU}(d)\rarw\calO_{X_1}(d)\rarw 0,$$ 
where for any sheaf $\cal{F}$ we denote 
$\calF(d)=\calF \otimes \calO_{\CPTU}(d)$. 
Passing to cohomologies we obtain: 
\[ \begin{CD} 
0 @>>> H^0(\calJ_{X_1}(d)) @>>> H^0(\calO_{\CPTU}(d)) 
@>{R_{X_1}}>> H^0(\calO_{X_1}(d)) @>>> H^1(\calJ_{X_1}(d)) @>>> \ldots 
\end{CD} \] 
where the map ${R_{X_1}}$ is induced by the restriction 
$\calR_{X_1}: \calO_{\CPTU} \rarw \calO_{X_1}$. 
Since $H^1(\CPTU,\calJ_{X_1}(d)) \cong 
H^1(\wtldS_p, \calO_{\wtldS_p}(D-mE))$ and the latter vanishes by assumption 
we see that the map $R_{X_1}$ is surjective.

The choice of the affine chart $\CC^2 \subset \CPTU$ induces an isomorphism 
$i_1:\calP(d) \rarw H^0(\calO_{\CPTU}(d))$, where $\calP(d)$ denotes 
the space of polynomials in $(x,y)$ of degree not more than $d$. 
Similarly, we obtain an isomorphism $i_1':\oplus_{i=1}^n J^{(m_i)}_{p_i} 
\oplus J^{(m)}_{p_0} \rarw H^0(\calO_{X_1})$.
Denoting by 
$$\wtld{R}_1 : \calP(d) \rarw \bigoplus_{i=1}^n J^{(m_i)}_{p_i} 
\bigoplus J^{(m)}_{p_0}$$ the linear map 
$$\wtld{R}_1(F) = 
\left(j^{(m_1)}_{p_1}(F), \ldots, j^{(m_n)}_{p_n}(F), j^{(m)}_p(F)\right),$$ 
we obtain the following commutative diagram: 

\[
\begin{CD}
\calP(d) @>{\wtld{R}_1}>> \bigoplus_{i=1}^n J^{(m_i)}_{p_i} 
\bigoplus J^{(m)}_{p_0} \\
@V{i_1}VV  @VV{i_1'}V \\ 
H^0(\calO_{\CPTU}(d))  @>{R_{X_1}}>>  H^0(\calO_{X_1}(d))
\end{CD}
\]

As $R_{X_1}$ is surjective so is $\wtld{R}_1$.
But $\wtld{R}_1^{-1}(\oplus_{i=1}^n J^{(m_i)}_{p_i})= \calP(\Delta_1)
\subset \calP(d)$ and ${\wtld{R}_1}|_{\calP(\Delta_1)} = R_1$. This implies 
that $R_1$ is indeed surjective.

The case of $R_2$ is easier. Replacing $d$ by $m'$, $X_1$ by the subscheme 
$X_2\subset \CPTU$ defined by the ideal sheaf 
$\calJ_{X_2}=\prod_{j=1}^k \calJ^{\al_j}_{q^0_j}$, and $R_{X_1}$ by 
the the restriction map $R_{X_2}$, we obtain the commutative diagram:  

\[
\begin{CD}
\calP(\Delta_2) @>{R_2}>> \bigoplus_{j=1}^k J^{(\al_i)}_{q^0_i}  \\
@V{i_2}VV  @VV{i_2'}V \\ 
H^0(\calO_{\CPTU}(m'))  @>{R_{X_2}}>>  H^0(\calO_{X_2}(m'))
\end{CD}
\]
where $i_2$ and $i_2'$ are obvious isomorphisms induced by the choice of the 
affine chart $\CC^2 \subset \CPTU$.
The vanishing of $H^1(V,\calO_{V}(C)) \cong H^1(\CPTU,\calJ_{X_2}(m'))$ 
implies, as before, the surjectivity of $R_{X_2}$ and consequently that of 
$R_2$.

To conclude the proof of the lemma, consider the smooth family of linear maps 
$$R^{(t)}:\calP(\Delta) \rarw  \bigoplus_{i=1}^n J^{(m_i)}_{p_i} 
\;\;\bigoplus \;\;
\bigoplus_{j=1}^k J^{(\al_j)}_{q^0_j},$$ defined by:\footnote{As before,
using the convention that $t^0 \equiv 1$ the family $R^{(t)}$ extends also 
for $t=0$.}
$$R^{(t)}\big(\sum_{(i,j)\in \Delta} A_{ij} x^i y^j\big) = 
\left(R_1\big(\sum_{(i,j) \in \Delta} A_{ij} x^i y^j t^{\nu_1(i,j)}\big), 
\quad R_2\big(\sum_{(i,j) \in \Delta} A_{ij} x^i y^j t^{\nu_2(i,j)}\big) 
\right).$$ 

Substituting $t=0$ we have, under the direct sum decomposition 
$\calP(\Delta)=\calP(\Delta_1) \oplus \calP(\Delta_2)$, that 
$R^{(0)}=R_1 \oplus R_2$, hence $R^{(0)}$ is surjective.
Since the family $R^{(t)}$ depends smoothly on $t$ we conclude that $R^{(t)}$ 
remains surjective for $t>0$ small enough. By the (linear) implicit function 
theorem there exists a smooth deformation $\{A_{ij}(t)\}_{0\leq t \leq \eps}$ 
of $a_{ij}$, such that 
$R^{(t)}(\sum_{(i,j)\in \Delta} A_{ij}(t)x^i y^j)=0$. This means 
that $F_{1,t}(x,y)$ vanishes to order $m_i$ at $p_i$ for every 
$1\leq i \leq n$ and $F_{2,t}(x,y)$ vanishes to order $\al_j$ at $q_j^0$ 
for every $1\leq j \leq k$. 
This concludes the proof of the lemma and thus of the whole proposition.
\ \Qed

\subsection{Passing from specific points to very general} \cntrsb

In what follows we shall detect several useful ample (resp. nef) vectors
$(d;\al_1, \ldots,\al_k)$ by choosing $k$ points $q_1,\ldots,q_k\in
\CPTU$ to lie in a very specific convenient position which is not
generic. The following lemma shows that this vectors remain ample (resp.
nef) also for a very general choice of the points $q_1,\ldots,q_k$.

\begin{lem}
\label{lem-very_general}

Let $F$ be a divisor on a simple rational surface $S$, and $q^{(0)}_1,
\ldots q^{(0)}_k \in S \setminus (\Sig^S \cup Supp\,F)$ distinct
points. Let $\pi_0:\wtldS_0\rarw S$ be the blow-up of $S$ at
$q^0_1,\ldots q^0_k$ with exceptional divisors
$E^0_i=\pi^{-1}_0 (q^0_i) \, i=1,\ldots,k$.

Suppose that for some $f_1,\ldots,f_k \geq 0$ the divisor
$\pi^*_0F -\sum_{j=1}^k f_j E^0_j$ is ample (resp. nef). Then, for a
very general choice of points $q_1,\ldots,q_k \in S \setminus
(\Sig^S\cup Supp\,F)$ the divisor $$\pi^*F-\sum_{j=1}^k f_j E_j$$ is
ample (resp. nef) on the blow-up $\pi:\wtldS \rarw S$ of $S$ at
$q_1,\ldots, q_k$ with exceptional divisors $E_j=\pi^{-1}(q_j) 
\quad j=1,\ldots,k$.

\end{lem}

\pf 
The idea of the proof is very simple. 
Since $\wtld{F}_0=\pi^*_0F -\sum_{j=1}^k f_j E^0_j$ is assumed to be nef on 
the blow-up of $S$ at $q^0_1,\ldots q^0_k$, all the divisor classes which 
intersect $\wtld{F}_0$ negatively do not admit any effective representatives. 
Now, the point is that if a divisor class on the blow-up of $S$ at specific 
points has no effective representatives then the same will continue to hold 
also on the blow-up at generic points. The lemma now follows because 
$\mbox{Pic}(\wtldS)$ is countable. Let us give now the precise details. 

We prove the lemma for the ``nef'' case, the ``ample'' being very 
similar. Consider the following subset of 
$\mbox{Pic}(S)\times \ZZ_{\geq 0}^k$: 
$$\calB=\left\{(A,a_1,\ldots a_k) \in \mbox{Pic}(S)\times \ZZ_{\geq 0}^k 
\bigg| F\cdot A - \sum_{j=1}^k f_j a_j < 0 \right\}.$$
Notice that $\calB$ is a countable set.
We claim that for every $B=(A,a_1,\ldots a_k)\in \calB$ there exists a 
non-empty Zariski-open subset 
$\calU_B \subset \calC_k(S\stmin (\Sig^S \cup Supp\,F))$ such that for 
every $(q_1, \ldots, q_k) \in \calU_B$, the surface $\wtldS$ obtained by 
blowing-up $S$,\ $\pi:\wtldS \rarw S$, at $q_1, \ldots, q_k$ does not 
admit any effective divisor in the class $\pi^*[A]-\sum_{j=1}^k a_j [E_j]$.
Once this is proved, we take $\calV=\cap_{B\in \calB} \calU_B$. 
Obviously $\calV$ is a very general subset of  
$\calC_k(S\stmin (\Sig^S \cup Supp\,F))$ having the needed properties.

Let us prove the existence of the Zariski-open sets $\calU_B$ claimed 
above. For this end put $\calC=\calC_k(S\stmin (\Sig^S \cup Supp\,F)), \quad 
X=\calC \times S$, and denote by $pr:X \rarw \calC$ the obvious projection. 

Consider the subvarieties $Y_j\subset X, \quad j=1,\ldots,k$, defined
by $$Y_j=\{ ((x_1, \ldots, x_k),x)\mid x=x_j \}.$$ The $Y_j$'s are smooth 
disjoint subvarieties of $X$ each of which is mapped by $pr$ 
isomorphically onto $\calC$. Let $\Th: \wtld{X} \rarw X$ be the blow-up of 
$X$ along $Y=\cup_{j=1}^k Y_j$ and write $\wtld{E}_j=\Th^{-1}(Y_j)$ for the 
exceptional divisors.  

Given $B=(A,a_1,\ldots a_k)\in \calB$, we denote by $\calL$ the line bundle 
$$\calL=\calO_{\wtldX}(\wtld{A}-\sum_{j=1}^k a_j \wtld{E}_j) \in 
\mbox{Pic}(\wtldX),$$ where $\wtld{A}\in \mbox{Div}(\wtldX)$ is the divisor 
$\Th^*(\calC \times A)$. Finally, for every $\undl{q} \in \calC$ we 
write $\calL_{\undl{q}}$ for the restriction of $\calL$ to the surface 
$\wtldS_{\undl{q}} = \Th^{-1} pr^{-1}(\undl{q})$.

Let $\undl{q}=(q_1, \ldots, q_k)\in \calC$. It is easy to see that the 
map $\pi_{\undl{q}}:\wtld{S}_{\undl{q}} \rarw S$ defined by the composition 
$\wtldS_{\undl{q}} \stackrel{\Th}{\lrarw} X 
\stackrel{pr_{_S}}{\lrarw} S$ 
is just the blow-up of $S$ at $q_1, \ldots, q_k$, and that 
$$\calL_{\undl{q}}=\calO_{\wtld{S}_{\undl{q}}}(\pi^*_{\undl{q}} A-\sum_{j=1}^k 
a_j E_j).$$

Now, for $\undl{q}^0=(q_1^0, \ldots, q_k^0)$ we know that 
$\dim H^0(\wtld{S}_{\undl{q}^0},\calL_{\undl{q}^0})=0$ because 
$\pi_0^*F-\sum_{j=1}^k f_j E_j$ is nef. It follows from the 
semicontinuity theorem (see~\cite{Ha-AG}) that there exits a Zariski-open 
neighborhood of $\undl{q}^0$, \ $\calU_B \subset \calC$ such that for every 
$\undl{q} \in \calU_B$, \ $H^0(\wtld{S}_{\undl{q}},\calL_{\undl{q}})=0$. 
\ \Qed 

\subsection{Proof of the gluing Theorem} \cntrsb 
\label{subsect-prf_glue_thm}

Now we are in position to prove Theorem~\ref{thm-glue2}. 
\pf We divide the proof into three steps. 
In the first step we prove that the resulting divisor 
$\wtldD=\pi^*D-\sum_{j=1}^k \al_j E_j$ is nef provided that 
$D_p=\pi_p^*D-mE$ and 
$v=(m;\al_1,\ldots, \al_k)$ are nef. In the second step we prove that the 
theorem holds under the assumption that both $D_p$ and $v$ are ample. 
Finally, in the third step we prove the theorem in its full generality by 
reducing the problem to the first two steps.

{\em Step 1.} Assuming that $D_p=\pi_p^*D-mE$ and $v=(m;\al_1,\ldots,\al_k)$ 
are nef we prove that $\wtldD=\pi^*D-\sum_{j=1}^k E_j$ is nef.

We claim that there exists $N_0>0$ and a divisor $A$ on $S$ such that for 
every $N>0$ there exists a very general subset 
$G_N\subset \calC_k (S \stmin (\Sig^S\cup Supp\, D))$ such that for every 
$(q_1,\ldots, q_k) \in G_N$ the divisor 
$$D_N=\pi^*A+(N+N_0)\pi^*D - N \sum_{j=1}^k \al_j E_j$$ is nef on the blow-up 
$\pi:\wtldS \rarw S$ of $S$ at $q_1, \ldots, q_k$.

Once this is proved step 1 of the proof will be concluded as follows: put 
$G=\cap_{N=1}^{\infty} G_N$. 
Clearly $G \subset \calC_k (S \stmin (\Sig^S\cup Supp\, D))$ 
is a very general subset. Let $(q_1, \ldots, q_k) \in G$ and let 
$C \subset \wtldS$ be a curve. Since $D_N$ is nef we have 
$$0\leq D_N \cdot C =A\cdot C +(N+N_0)D\cdot C -N\sum_{j=1}^k C\cdot E_j.$$
Dividing by $N$ and letting $N\rarw \infty$ we obtain that 
$$\wtldD\cdot C = (\pi^*D-\sum_{j=1}^k \al_j E_j)\cdot C \geq 0.$$ 
As $v$ is nef, we have $m^2\geq \sum_{j=1}^k \al_j^2$ and so 
$$\wtldD\cdot \wtldD = D \cdot D - \sum_{j=1}^k \al_j^2 
\geq D\cdot D-m^2 = D_p \cdot D_p \geq 0,$$ the latter inequality following 
from the nefness of $D_p$. Thus $\wtldD$ is nef.

Let us prove now the existence of $N_0,A,G_N$ claimed above. For the divisor 
$A$ we choose any divisor on $S$ such that $\pi_p^*A-E$ is ample on 
$\wtldS_p$. The nefness of $v$ means by definition that there exit $k$ 
distinct points $q^0_1, \ldots, q^0_k \in \CPTU$ such that the 
divisor $mL^V-\sum_{j=1}^k \al_j E^V_j$ is nef on $V$ -- the blow-up 
of $\CPTU$ at $q^0_1, \ldots, q^0_k$. For $N_0>0$ we choose any 
integer for which $B=N_0 L^V - \sum_{j=1}^k E^V_j$ is ample on $V$.

For every $N>0$ define $L_N'\in 
\mbox{Div}(\wtldS_p),\,\, L_N''\in \mbox{Div}(V)$ to be:
$$L_N'=(\pi_p^*A-E)+(N+N_0)D_p=
\pi_p^*A+(N+N_0)\pi_p^*D-(mN+mN_0+1)E,$$
$$L_N''=B+N(mL^V-\sum_{j=1}^k \al_jE_J^V)=
(Nm+ N_0)L^V-N\sum_{j=1}^k E^V_j.$$ 
It easily follows from our assumptions on $D_p$ and on $v$ that $L_N',L_N''$ 
are ample for every $N>0$.
Choose an intger $r_N>0$ for which the following two conditions are 
satisfied: 
\begin{enumerate}
\item[1)] $r_N L_N'$ and $r_N L_N''$ are very ample. 
\item[2)] $H^1(\wtldS_p, \calO_{\wtldS_p}(r_N L_N'))=0$, and 
$H^1(V, \calO_V(r_N L_N''))=0$. 
\end{enumerate}

Choose irreducible curves $\wtldC_N' \in |r_N L_N'|$ and 
$\wtldC_N'' \in |r_N L_N''|$ and put $C_N'=\pi_p(\wtldC_N') \subset S$.
By Proposition~\ref{prop-glue} there exist 
$q_1, \ldots, q_k \in S\stmin (\Sig^S \cup C_N')$ such that the surface 
$\wtldS$ obtained by blowing-up $\pi:\wtldS \rarw S$ at $q_1, \ldots, q_k$ 
admits an irreducible curve $C_N$ in the linear system
$$C_N\in \left| r_N 
\bigg(\pi^*A+(N+N_0)\pi^*D - N \sum_{j=1}^k \al_j E_j\bigg) 
\right| = \left| r_N D_N \right|.$$ 
Noting that 
$$D_N\cdot D_N \geq N^2(D \cdot D - \sum_{j=1}^k \al_j^2) 
\geq N^2(D\cdot D -m^2)=N^2 D_p\cdot D_p \geq 0,$$ we conclude that $D_N$ 
intersects every curve non-negatively and so it is is nef on $\wtldS$. 
By Lemma~\ref{lem-very_general} we may assume that 
$(q_1, \ldots,q_k)$ vary in some very general subset 
$G_N\subset \calC (S \stmin (\Sig^S\cup Supp\, D))$. This completes the 
proof of step 1.

{\em Step 2.} Assuming $D_p$ and $v$ are both ample we prove that 
$\wtldD=\pi^*D-\sum_{j=1}^k \al_j E_j$ is ample.

Here it is more convenient to work with $\QQ$-divisors. First note that 
step 1 remains true if we take $m$ and $\al_j$ to be rational numbers.
It follows from Seshadri's criterion for ampleness (see\cite{Ha-Ample}) 
that there exists a positive rational number $\eps$ such that 
$\pi^*D-(1+\eps)mE$ is ample. 
Since $( (1+\eps)m; (1+\eps)\al_1, \ldots, (1+\eps)\al_k )$ is ample too we 
have from step 1 that 
$\wtldD_{\eps}=\pi^*D-(1+\eps)\sum_{j=1}^k \al_j E_j$ is nef.

Let us prove that $\wtldD$ is ample by applying Nakai-Moishezon criterion. 
Indeed, let $\wtldC \subset \wtldS$ be a curve.
If $\wtldC=E_j$ is one of the standard exceptional divisors then 
$\wtldD \cdot \wtldC = \al_j > 0$ (recall that $\al_j >0$ because we assume 
that $v$ is ample). Otherwise, let $C=\pi(\wtldC)$.
If $\wtldC$ does not pass through any of the exceptional divisors then 
$\wtldD\cdot \wtldC= D\cdot C >0$ because $D$ itself is ample for $D_p$ is.
Suppose now that there exits a $j_0$ such that $\wtldC\cdot E_{j_0}>0$. 
In this case $\wtldD \cdot \wtld C = 
D\cdot C - \sum_{j=1}^k \al_j \wtldC\cdot E_j > 
\wtldD_{\eps} \cdot \wtldC \geq 0$. Finally note that 
$\wtldD\cdot \wtldD > \wtldD_{\eps}\cdot \wtldD_{\eps} \geq 0$.

{\em Step 3.} Consider the general case. 
The case of $D_p$ nef has been treated in step 1 so we may assume that 
$D_p$ is ample and $v$ is nef. Similarly to step 2 we choose 
a positive rational number $\eps$ such that both $\pi_p^*D-(1+\eps)E$ and 
$((1+\eps)m; \al_1, \ldots, \al_k)$ are ample. By step 2 we have that 
$\pi_p^*D-\sum_{j=1}^k E_j$ is ample.
\ \Qed
\ \\ 

Theorem~\ref{thm-glue1} follows immediately from Theorem~\ref{thm-glue2} by 
taking $D=dL^S-\sum_{i=1}^n m_i E^S_i$ for a suitable simple rational 
surface $S$. 
\ \Qed 

\section{Asymptotics on the remainders of $\Oone$} \cntrs
\label{sect-Asymptotics}

In order to obtain estimates on $\calR_N(\Oone)$ we shall extensively
use Theorem~\ref{thm-glue1} in combination with the
{\em Cremona action}.
The point is, that the Cremona group acts on the set of ample (resp.
nef) vectors. Let us briefly summerize the needed facts about the
Cremona action. We refer the reader to~\cite{Do-Or} for more details.

\subsection{The Cremona action on the ample cone} \cntrsb
\label{subsect-Cremona}

Denote by $(H_k,\langle \, , \, \rangle) \quad (k\geq 3)$ the 
hyperbolic lattice
$H_k=\ZZ l \oplus\ZZ e_1 \oplus \ldots \oplus \ZZ e_k$
with the bilinear form $\langle\, ,\, \rangle$ defined by 
$\langle l,l\rangle=1, \quad \langle l, e_j\rangle=0, \quad 
\langle e_i,e_j\rangle=-\dlt_{ij}$. Consider the subgroup 
$Cr_k \subset Aut(H_k,\langle \, , \,\rangle)$,
generated by: 
\begin{enumerate}
\item[1)] The symmetric group 
$S_k \hookrightarrow Aut(H_k,\langle\, ,\,\rangle)$ acting
on the last $k$ components. 
\item[2)] The reflection 
$R_{123}:(H_k,\langle \, ,\,\rangle) \rarw (H_k,\langle \, ,\,\rangle)$ 
defined by $R_{123}(\eta)=\eta + \langle \eta, r_{123}\rangle r_{123}$, 
where $r_{123}=l-e_1-e_2-e_3$. 
\end{enumerate}
The group $Cr_k$ is called the {\em Cremona group}.

It is easily seen that the reflections
$R_{ijk}(\eta)=\eta + \langle \eta, r_{ijk}\rangle r_{ijk}$, where
$r_{ijk}=l-e_i-e_j-e_k$, belong to $Cr_k$. Let us mention one more useful 
transformation which we denote by $SR$. The transformation
$SR$ takes a vector $v=(d;m_1, \ldots, m_k) \in H_k$ and sorts it. In other 
words $SR(v)=(d;m_{\tau(1)}, \ldots, m_{\tau(k)})$, where $\tau$ is a 
permutation of $\{1,\ldots, k\}$ for which 
$m_{\tau(1)} \geq  \ldots \geq m_{\tau(k)}$. It is obvious that for every 
vector $v\in H_k$ there exists $\sig \in Cr_k$ such that $SR(v)=\sig(v)$.

Given a simple rational surface obtained by blowing up $\Th:V\rarw
\CPTU$ of $p_1,\ldots,p_n \in \CPTU$, there is an isomorphism of
lattices $m_{\Th}:(\mbox{Pic}(V),\,\cdot \,) \rarw 
(H_k,\langle \, ,\, \rangle)$, 
where $\cdot$ stands for the intersection form on $\mbox{Pic}(V)$. 
The isomorphism $m_{\Th}$ sends $L^V$ to $l$ and $E^V_i$ to $e_i$.

To deduce that $Cr_k$ acts on the set of ample (resp. nef) vectors we
need the following lemma which essentially appears in~\cite{Do-Or}.

\begin{lem}
\label{lem-marking}

Let $V$ be a simple rational surface obtained by blowing-up $\Th:V\rarw
\CPTU$ points $p_1,\ldots,p_n \in \CPTU$ in general position.
Then for every $\sig \in Cr_k$ there exists a simple rational surface
$V_{\sig}$ obtained by blowing-up $\Th_{\sig} : V_{\sig} \rarw \CPTU$
points $q_1,\ldots q_k$ in general position and a biholomorphism
$f_{\sig}:V_{\sig} \rarw V$ making the following diagram commutative:
\[
\begin{CD}
\mbox{\em Pic}(V) @>{f^*_{\sig}}>> \mbox{\em Pic}(V_{\sig}) \\
@V{m_{\Th}}VV  @V{m_{\Th_{\sig}}}VV \\
H_k @>{\sig}>>  H_k
\end{CD}
\]

\end{lem}

Combining this with Lemma~\ref{lem-very_general} we immediately obtain
the following

\begin{lem}
\label{lem-Cremona}

When $k\geq 3$ the group $Cr_k$ acts on the set of ample (resp. nef)
vectors viewed as a subset of $H_k$.

\end{lem}

\noindent 
{\bf Remark.} From Lemma~\ref{lem-very_general} it follows that there 
exists (at least) one simple rational surface $S$, obtained by
blowing-up $k$ distinct points in $\CPTU$, $\Th:\wtldS\rarw \CPTU$ such
that $\calL \in \mbox{Pic}(S)$ is ample (resp. nef) iff $m_{\Th}(\calL)$ is 
ample (resp. nef). 
Hence, the set of ample (resp. nef) vectors is closed under addition and 
multiplication by positive (resp. non-negative) integers. 
Henceforth we shall denote by $\calK_k \subset H_k \otimes \RR$ 
(resp. $\calKbar_k$) the cone generated by all ample (resp. nef) vectors.

\subsection{An algorithmic procedure for detecting ample classes} \cntrsb 
\label{subsect-Alg}

Given two vectors $v_1=(d;m_1, \ldots, m_n) \in H_n$ and
$v_2=(\dlt; \al_1, \ldots, \al_k) \in H_k$ with $\dlt=m_i$ for some
$1\leq i\leq n$, define a new vector $v_1 \dzi v_2 \in H_{n+k-1}$ by setting 
$$v_1 \dzi v_2=(d; m_1,\ldots, m_{i-1}, \al_1, \ldots, \al_k, m_{i+1}, 
\ldots, m_n).$$
Theorem~\ref{thm-glue1} asserts that if
$v_1$ is ample (resp. nef) and $v_2$ is nef, then $v_1 \dzi v_2$
is ample (resp. nef).

Given a vector $v_0 \in H_N$ the ampleness of which we want to prove
we try to find a decomposition $v_0=v_1 \#_{_{i_1}} u_1$ where $u_1 \in
H_{k_1}$ is known to be nef and $v_1\in H_{n_1},\,\,(k_1+n_1 - 1 = N)$.
If $v_1$ turns to be ample then we are done in view of
Theorem~\ref{thm-glue1}. To check the ampleness of $v_1$ we first 
"simplify" it by applying to it Cremona transformations. For example, we
may try, using Cremona transformations to reduce the degree of $v_1$
(by the degree of $v=(d;\mu_1, \ldots, \mu_k)$ we mean $\deg\,v=d$).
Let $v_1'$ be a simpler vector in the same orbit of $v_1$ under
the action of $Cr_{n_1}$ (e.g. $v_1'$ having minimal degree in the orbit, 
or having some other convenient feature). By Lemma~\ref{lem-Cremona} $v_1$ is
ample iff $v_1'$ is. Now we apply the whole process to $v_1'$ and so
on. In this way we obtain a sequence of vectors
$v_1,u_1,v_1', \ldots, v_r,u_r,v_r'$ where $v_j'$ is a Cremona
simplification of $v_j \in H_{n_j}$, $u_j \in H_{k_j}$ is a nef vector
and $v_j'=v_{j+1} \#_{_{i_{j+1}}} u_{j+1}$ for some $i_{j+1}$.

Note that at each stage the number of points decreases, namely
$n_{j+1}<n_j$ provided that $k_j > 1$. The process ends successfully as
soon as we are able to prove that $v_r$ is ample for some $r$.
We remark that if one of the $v_j$ turns out not to be ample 
then process fails to give any information because the converse of
Theorem~\ref{thm-glue1} is not true. However, we may attempt to find other 
decomposition sequences $v_1,u_1,v_1', \ldots,v_r,u_r,v_r'$ 
(see Section~\ref{sect-Remarks}).

The same procedure can be applied for proving nefness of a
vector $v_0$ by requiring that $v_r$ is nef instead of ample.
In the next subsection we shall apply this process in order to prove 
Theorem~\ref{thm-asymp1} and Corollary~\ref{cor-coef2}.

In order to make the preceding procedure applicable we must first endow 
ourselves with an initial large enough collection of ample and nef vectors
which will play the role of the $u_j$'s and of $v_r$. To simplify notations 
let us agree that 
$(d;\al_1^{\times r_1}, \ldots \al_k^{\times r_k})$ stands for 
$$(d;\underbrace{\al_1, \ldots, \al_1}_{\mbox{\tiny $r_1$ times}}, 
\ldots \ldots ,\underbrace{\al_k, \ldots, \al_k}_{\mbox{\tiny $r_k$ times}})
\in H_N,$$ where $N=\sum_{j=1}^k r_j$.

The next lemma provides a modest initial collection of ample and nef vectors 
which is sufficient for our purposes.

\begin{lem}
\label{lem-vectors}

The following vectors are nef (resp. ample) on a very general rational
surface: 
\begin{enumerate}
\item[1)] $(d; 1^{\times r})$, where $d^2 \geq r$ (resp. $d^2 > r$). 
\item[2)] $(d; m_1, m_2, 1^{\times r})$ where $d\geq m_1+m_2$ and 
$d^2 \geq m_1^2+m_2^2+r$.
\end{enumerate}


\end{lem}

\noindent
{\bf Remark.} The ``ample'' case of statement 1 above has been proved by Xu 
in~\cite{Xu-Divisors} and by K\"{u}chle in~\cite{Ku}. Below however, we 
present an alternative proof suggested by Ilya Tyomkin.

\pf 
Notice first that in view of Lemma~\ref{lem-very_general} it
is enough to prove that the above vectors are nef (resp. ample) on a specific 
simple rational surface.\\ \ \\
1) Consider first the case $d^2>r$. In~\cite{Nag} (consult also~\cite{Sh-Ty}) 
Nagata proved that if $N$ is a square, then for
generic points $p_1, \ldots, p_N\in \CPTU$ and for every irreducible
curve $C\subset \CPTU$ the following strict inequality holds:
$$\mbox{deg}(C) > \frac{\sum_{j=1}^N \mbox{mult}_{p_j}(C)}{\sqrt{N}}.
\eqno(1)$$
Let $V_r$ be the blow-up of $\CPTU$ at $r$ generic points and denote by
$\Th:\wtld{V}_r \rarw V_r$ the blow-up of $V_r$ at $d^2-r$ generic
points. Thus $\wtld{V}_r$ is the blow-up of $\CPTU$ at $N=d^2$
generic points and it follows from inequality (1) that the divisor
$\wtld{D}=\Th^*(dL^{V_r}-\sum_{j=1}^r E^{V_r}_j) - \sum_{j=r+1}^{d^2} E_j$ 
intersects every curve positively. This immediately implies that 
$D=dL^{V_r}-\sum_{j=1}^r E^{V_r}_j$ intersects any curve in $V_r$ positively. 
As $D\cdot D > 0$ the statement follows from Nakai-Moishezon criterion 
(see~\cite{Ha-Ample}).

The proof for the nef case ($d^2 \geq r$) is much easier. Indeed, let
$C\subset \CPTU$ be an irreducible curve of degree $d$, and let $p_1,
\ldots, p_r$ be distinct points on $C$ at which $C$ is smooth. Let $V$
be the blow-up of $\CPTU$ at $p_1, \ldots, p_r$ and let $D$ be the
proper transform of $C$ in $V$, $D\in | d L^V-\sum_{j=1}^r E_j|$.
As $D$ is an irreducible curve of non-negative self intersection the
vector $(d; 1^{\times r})$ corresponding to the divisor class of $D$ is nef on
$V$. 
\ \\ \ \\
2) Set $D=dL - m_1 E_1 - m_2 E_2$. Consider the linear system $|D|$ on
$V_2$ -- the blow-up of $\CPTU$ at 2 points.
As $D=m_1(L-E_1) + m_2(L-E_2) + (d-m_1-m_2)L$ it is easy to see that
$|D|$ is not empty and has no base-points, hence by Bertini theorem
there exists an irreducible (smooth) curve $C\in |D|$. Choose $r$
distinct points $p_1, \ldots, p_r \in C \setminus (E_1 \cup E_2)$ and
let$\wtld{V}$ be the blow-up of $V$ at $p_1, \ldots, p_r$. Finally
denote by $\wtld{C}$ be the proper transform of $C$ in $\wtld{V}$.

We have $\wtld{C} \in |dL- m_1E_1 -m_2E_2 - \sum_{j=3}^{r+2} E_j|$.
As $\wtld{C}$ is irreducible and $\wtld{C} \cdot \wtld{C} \geq 0$, the
vector $(d;m_1,m_2,1^{\times r})$ is nef.
\ \Qed
\ \\ \\
{\bf Remark.} Note that the cones $\calK_n$ (resp. $\calKbar_n$) can be
explicitly computed when $n<9$ (see~\cite{Dmz},~\cite{F-M}), and so 
can be joined to the initial collection of ample and nef vectors to be
applied in the framework of the process mentioned above.

\subsection{Proofs of Theorem~\ref{thm-asymp1} and Corollary~\ref{cor-coef2}} 
\cntrsb 
\label{subsect-prfs1} 

We start with the proof of Theorem~\ref{thm-asymp1}. 
\pf 1) Let $N=a^2l^2+2l$ and $v_0=(a^2l+1;a^{\times N})$. 
As $\langle v_0, v_0 \rangle = 1$, nefness of $v_0$ will 
give the needed estimate for $\calR_N(\Oone)$. 

The decomposition $N=(al-1)^2 + n$, where $n=2al+2l-1$, leads us 
to $v_0=v_1 \#_{_1} u_1$ where 
$$v_1=\Big( a^2l+1; a(al-1), a^{\times n} \Big) \in H_{n+1}, \quad 
  u_1=a\Big(al-1; 1^{\times (al-1)^2}\Big)\in H_{(al-1)^2}.$$ 
By Lemma~\ref{lem-vectors} $u_1$ is nef, hence in view of 
Theorem~\ref{thm-glue1} we are reduced to proving that 
$v_1$ is nef. This turns out to be easy by using Cremona transformations.
Indeed let 
$$v_1'=R_{1,n-1,n} \circ R_{1,n-3,n-2} \circ \ldots \circ R_{123} (v_1),$$  
where $R_{ijk}\in Cr_{n+1}$ are defined in Section~\ref{subsect-Cremona}. 
A straight forward computation shows that 
$v_1'=(a+l; l-1, 1^{\times n-1},a)$. This vector is nef by 
Lemma~\ref{lem-vectors}, and therefore $v_1$ too. \\

\noindent 
2) Let $N=a^2l^2-2l$ and $v_0=(a^2l-1; a^{\times N})$. 
Again $\lbr v_0,v_0\rbr = 1$, hence in order to prove the needed estimate 
on $\calR_N(\Oone)$ we have to prove that $v_0$ is nef. 
Using the decomposition $N=(al-2)^2+n_1$, where $n_1=4al-2l-4$,  
we note that $v_0=v_1 \#_{_1} u_1$ where 
$$v_1=\Big( a^2l-1; a(al-2), a^{\times n_1} \Big) \in H_{n_1+1}, \quad
  u_1=a\Big( al-2; 1^{\times (al-2)^2}\Big) \in H_{(al-2)^2}.$$ 
By Lemma~\ref{lem-vectors} $u_1$ is nef. We are thus reduced to proving 
nefness of $v_1$. By applying similar Cremona transformations as in 1, we 
obtain that $v_1'=(a^2l-2al+l+1; (al-l-2)(a-1), (a-1)^{\times n_1})$ lies in 
the same orbit as $v_1$.

Let us apply now the same algorithm again on $v_1'$. For this, consider the 
decomposition $v_1'=v_2\#_{_1} u_2$, where 
$$v_2=\Big(a^2 l-2al+l+1; (al-l-1)(a-1), (a-1)^{\times 2al-1}\Big) 
\in H_{2al-1},$$
$$u_2=(a-1)\Big( al-l-1;  al-l-2, 1^{\times 2al-2l-2}\Big).$$ 
By Lemma~\ref{lem-vectors} $u_2$ is nef, thus we are reduced to proving that 
$v_2$ is nef. Using similar Cremona transformations as in 1 we 
obtain that $v_2'=\big(a+l-1; l-1, 1^{\times 2al-2}, a-1 \big)$ lies in the 
same orbit as $v_2$. But by Lemma~\ref{lem-vectors} $v_2'$ is nef. \\ 

3) Let $N=a^2 l^2+l$ and suppose that $a=2^k b$ with $k\geq 0$ and $b$ odd.
The assumption appearing in the statement of the Theorem is that $l>2b$. 
Note that we may assume that $l$ is odd, since when $l$ is even we have 
$N=(2a)^2(\frac{l}{2})^2+2\frac{l}{2}$ and this is already covered 
in 1 above. 

In order to prove the needed estimate on $\calR_N(\Oone)$ we have to show 
that the vector $v=\big(2a^2 l+1; 2a^{\times (a^2l^2+l)}\big)$ is nef.
Let us prove a slightly stronger statement, namely: 

{\em Claim.} The vector $v_0=\big(2a^2 l+1; 2a^{\times (a^2l^2+l)}, 1\big)$  
is nef. 

We argue by induction on $k$. Consider first the case $k=0$. 
We have $v_0=w\#_{_1}u$, where 
$$w=\Big(2a^2l+1; 2a(al-1), 2a^{\times 2al+l-1},1\Big),
\quad \mbox{and} \quad 
u=2a\Big(al-1; 1^{\times (al-1)^2}\Big).$$ 
The latter being nef, we are reduced 
to proving nefness of $w$. Applying suitable Cremona transformation to $w$, 
we obtain the vector 
$$w'=\Big(\frac{l+1}{2}+a; \frac{l-1}{2}-a, 1^{\times 2al+l}\Big).$$
Since $l>2b=2a$ we have that $\frac{l-1}{2}-a\geq 0$ and so $w'$ is nef by 
Lemma~\ref{lem-vectors}. This completes the basis of the induction. 

Let us turn now to the case $k>0$. We have 
$$v_0=\Big(\,\big(v_1\#_{_3}u_1\big)\#_{_2}u_1\Big)\#_{_1}u_1,$$ where 
$$v_1=\Big(2a^2l+1;(a^2l)^{\times 3},2a^{\times (\frac{a}{2}l)^2+l}, 1\Big), 
\quad 
u_1=2a\Big(\frac{a}{2} l; 1^{\times (\frac{a}{2}l)^2}\Big).$$ 
Again, $u_1$ is nef.
As for $v_1$, it lies in the same orbit under the Cremona action as the 
vector $v_1'=\big(a^2l+2; 2a^{\times (\frac{a}{2}l)^2+l}, 1^{\times 4}\big)$. 
Consider now the decomposition $v_1'=v_2 \# (2; 1^{\times 4})$, where 
$v_2=\big(a^2l+2; 2a^{\times (\frac{a}{2}l)^2+l},2\big)$ and $\#$ stands for 
gluing at the last coordinate of $v_2$. As $(2; 1^{\times 4})$ is nef, it is 
enough to prove that $v_2$ is nef. 
Putting $c=\frac{a}{2}=2^{k-1}b$, we have that 
$$v_2=2(2c^2l+1; 2c^{\times c^2 l^2+l},1).$$ By the induction hypothesis 
$v_2$ is nef. This completes the proof of the claim. The Theorem now follows 
easily. 
\ \Qed \ \\

Let us turn now to the proof of Corollary~\ref{cor-coef2}.
\pf 
Let $D=\pi^*\calO_{\CPTU}(d)-2\sum_{j=1}^N E_j$ and suppose that 
$D \cdot D \geq 0$. 

{\em Step 1.} Consider first the case $N=k^2+k$ for some $k$. 
By Theorem~\ref{thm-asymp1}-3 $$\calR_N(\Oone)\leq \frac{1}{(2k+1)^2}.$$ 
Since $D\cdot D =1$, this implies that $D$ is nef. 

{\em Step 2.} Consider the general case.
The condition $D\cdot D \geq 0$ reads $d^2\geq 4N$. We may assume that $d$ is 
odd, for the case of $d$ even is precisely the contents of Xu's theorem from 
Section~\ref{sect-Intro} (see~\cite{Xu-Divisors}. 
Writing $d=2k+1$, the condition $d^2\geq 4N$ gives $k^2+k >N$. 
By step 1, $\pi^*\calO_{\CPTU}(d)-2\sum_{j=1}^{k^2+k}E_j$ is nef, hence also 
$\pi^*\calO_{\CPTU}(d)-2\sum_{j=1}^N E_j$.
\ \Qed
\ \\ \\
{\bf Remark.} More careful considerations, in the spirit of the proof of 
Theorem~\ref{thm-asymp1} actually show that {\em when $d>5$, the divisor 
$D=\pi^*\calO_{\CPTU}(d)-2\sum_{j=1}^N E_j$ is ample iff $D\cdot D>0$.}

To prove this one has to sharpen first the second statement of 
Lemma~\ref{lem-vectors} and prove that $(d; m_1, m_2, 1^{\times r})$ is ample 
when $d >  m_1+m_2$ and $d^2 > m_1^2+m_2^2+r$. This can be done by 
similar, though more delicate, arguments to those used to prove nefness of 
these vectors. Then, using the ``ample+nef $\Rightarrow$ ample'' case of 
Theorem~\ref{thm-glue1} one deduces as in the proof of 
Theorem~\ref{thm-asymp1} that the divisor $D$ is ample for 
$N=k^2+k$, when $k>2$. The case of general $N$ can be easily reduced to 
$N=k^2+k$ as in the preceding proof.

\subsection{Proof of Corollary~\ref{cor-Nag}} \cntrsb
\label{subsect-prfs2}

\pf Let $N=N_1 N_2$. Nagata's conjecture for $N$ is equivalent to 
the nefness of vector $v=(d; m^{\times N})$ for every $d,m>0$ which satisfy 
$d^2-Nm^2>0$. 

Let $d,m$ be two such numbers. Choose a positive rational number $x$ such 
that $d^2 > x^2 N_2 > m^2 N$. The assumption of Nagata's conjecture for 
$N_1$ and $N_2$ implies that the vectors
$u=(x;m^{N_1}) \in \QQ^{N_1+1}$ and 
$w=(d; x^{\times N_2}) \in \QQ^{N_2+1}$ are nef. 

We have 
$v=\big( \ldots ((w \#_{_{N_2}}u) \#_{_{N_2-1}}u) \ldots \big)\#_{_{1}}u$.
Observing that Theorem~\ref{thm-glue2} remains valid also for vectors of 
rational numbers, we conclude that $v$ is also nef.
\ \Qed

\section{A conjecture relating continued fractions and remainders of
a line bundle} \cntrs
\label{sect-Conjecture}

The goal of this section is to propose a conjecture concerning estimates
on the values of the homogeneous remainders of $\Oone$, defined in
Section~\ref{subsect-Applic}. It turns out that all the cases appearing in
the statement of Theorem~\ref{thm-asymp1} are particular cases of this
conjecture.

Let us first recall some relevant facts from classical number theory.
Given a square-free number $N$, consider the following Diophantine
equation in the unknowns $d,m$ $$d^2-Nm^2=1.$$ This equation had been
attached-to the name {\em Pell's equation} in the ancient
literature and has been extensively studied by many
mathematicians in the 17'th and 18'th centuries including Leonard Euler
(see~\cite{Niv,Ir-Ro, VndP}). the classical result about the solutions
of this equation is that all solutions come from
{\em continued fractions expansions} of $\sqrt{N}$. Let us write
$\lbr a_0, a_1, \ldots, a_n \rbr$ for the continued fractions expansion
\[
a_0+
\cfrac{1}{a_1+\cfrac{1}{\ddots +\cfrac{1}{a_{n-1} +\cfrac{1}{a_n}}}}
\]

\noindent Similarly, we denote by $\lbr a_0, a_1, \ldots \rbr$ an infinite 
continued fractions expansion. It is not hard to see that the continued
fractions expansion of $\sqrt{N}$ must be of the following special
periodic form
$$\sqrt{N}=\lbr a_0,a_1,\ldots, a_{n-1}, 2a_0, a_1,\ldots, a_{n-1}, 2a_0,
\ldots \rbr,$$ hence we shall write from now on
$\sqrt{N}=\lbr a_0,\overline{a_1,\ldots,a_{n-1}, 2a_0} \rbr$ where the bar
denotes the periodic part. Moreover, it turns out that $a_i=a_{n-i}$ for
every $1\leq i \leq n-1$, (i.e. $(a_1,\ldots,a_{n-1})$ is a palindrome).

Define a rational number $\frac{d}{m}$ as follows: if $n$ is even put  
$$\frac{d}{m}=\lbr a_0, a_1,\ldots, a_{n-1} \rbr,$$ while for $n$ odd 
$$\frac{d}{m}=\lbr a_0, a_1,\ldots, a_{n-1}, 2a_0, a_1,\ldots, a_{n-1}\rbr.$$ 
It is well known that $(d,m)$ provides the minimal solution of Pell's 
equation, called the {\em fundamental solution}. 
Moreover, any other solution of Pell's
equation is obtained in a similar manner -- by truncating the infinite
continued fraction of $\sqrt{N}$ one term before the end of one of its
periods. More precisely, $(d,m)$ solves Pell's equation iff
$$\frac{d}{m}=\lbr a_0,\overline{a_1,\ldots, a_{n-1}, 2a_0}^{\, \times r}, 
a_1, \ldots, a_{n-1}\rbr , \,\,\, \mbox{where $r$ is odd if $n$ is odd}.$$
This formula means that the periodic 
part \ $a_1, \ldots a_{n-1}, 2a_0$ \ should be taken $r$ times and then once 
more without the last member $2a_0$. The number $r$ is allowed to be any 
non-negative integer in case $n$ is even, and $r$ must be odd if $n$ is odd. 
\ \\

Our conjecture is the following
\begin{cnj}
\label{cnj-CF}

Let $N>9$ be a square-free number and let $\frac{d}{m}$ be the fundamental 
solution of the corresponding Pell's equation. Then:\\ 
{\em 1)} The vector $(d;m^{\times N})$ is nef. \\
{\em 2)} $\calR_N(\Oone)\leq \frac{1}{d^2}$.

\end{cnj}

Our conjecture is much weaker than Nagata's, on the other hand it seems more 
accessible. Indeed, our methods provide a proof for the
conjecture in the following cases: \\ \ \\
1. {\em Consider the case that $\sqrt{N}$ has a 2-periodic continued
fractions expansion $\sqrt{N}=\lbr a_0,\overline{a_1,2a_0} \rbr$.} 
It is easy to see that this is the case iff $a_1|2a_0$ and $N=a_0^2 +
2\frac{a_0}{a_1}$. The solution of Pell's equation is 
$d=1_0a_1+1,\; m=a_1$. 
\begin{enumerate}
\item[A)] {\em Suppose that $a_1|a_0$}. Putting $a=a_1$ and 
$l=\frac{a_0}{a_1}$ we get $N=a^2l^2+2l$, and so by Theorem~\ref{thm-asymp1} 
our conjecture holds in this case. 
\item[B)] {\em Suppose that 
$a_1 \nmid a_0$ and $2^{k}a_0 > a_1^2$, 
where $k$ is the maximal integer for which $2^k | a_1$}. 
Since $a_1|2a_0$, $a_1$
must be even. Putting $a=\frac{a_1}{2}$ and $l=\frac{2a_0}{a_1}$
we obtain $N=a^2l^2+l$ and $l > \frac{a}{2^{k-2}}$. 
By Theorem~\ref{thm-asymp1} our conjecture holds. 
\end{enumerate} \ \\ 
2. {\em Consider $N$'s of the form $N=a^2l^2-2l$}. It is not hard to see
that $d=a^2l-1,m=a$ satisfy Pell's equation $d^2-Nm^2=1$. By
Theorem~\ref{thm-asymp1} we have $\calR_N(\Oone)\leq \frac{1}{d^2}$.
Therefore, if $(d',m')$ is the fundamental solution of Pell's equation
then $\calR_N(\Oone)\leq\frac{1}{d^2} \leq \frac{1}{{d'}^2}$, and so the
conjecture holds. Note that in this case the expansion of $\sqrt{N}$
will usually be longer than 2 
(example: $\sqrt{14}= \lbr 3,\overline{1,2,1,6} \rbr$).\footnote{It is not 
hard to see that if $N=a^2l^2-2l$ then $\sqrt{N}$ has 2-periodic expansion 
with minus signs.} \\ \ \\
3. Let us mention two other examples which do not fall into the above
categories. 
\begin{enumerate}
\item[A)] $N=19$. In this case 
$\sqrt{19}=\lbr 4,\overline{2,1,3,1,2,8} \rbr$. The fundamental solution is 
$d=170, m=39$. Thus, our conjecture suggests that 
$\calR_{19}(\Oone) \leq \frac{1}{{170}^2}$. 
\item [B)] $N=22$. In this case 
$\sqrt{22}=\lbr 4,\overline{1,2,4,2,1,8}\rbr$. The fundamental solution is 
$d=197, m=42$. Thus, our conjecture suggests that 
$\calR_{22}\leq \frac{1}{{197}^2}$. 
\end{enumerate}

Let us prove that the conjecture indeed holds in the cases 3.A and 3.B. 
We start with $N=19$.
By Lemma~\ref{lem-vectors} the vector 
$u=39(2; 1^{\times 4})=(78; 39^{\times 4})$ is nef. We have 
$$(170; 39^{\times 19})=
(\,(\,(170;78^{\times 3}, 39^{\times 7}) 
\#_{_{3}} u) \#_{_{2}} u) \#_{_{1}}u.$$ Thus, we are reduced to proving 
that $v=(170;78^{\times 3}, 39^{\times 7})$ is nef. To do this we apply the 
following Cremona transformations successively: 
\begin{enumerate}
\item[1)] Replace $v$ by $R_{123} (v)$. 
\item[2)] Sort the vector $v$, that is, replace $v$ by $SR(v)$, where the 
transformation $SR$ is the one defined in~\ref{subsect-Cremona}. 
\end{enumerate}
Applying this process enough times we finally arrive to the vector 
$(1; 0^{\times 10})$ which is nef. 

The case $N=22$ is similar. Here we use the decomposition 
$$(197; 42^{\times 22})=
(\,(\,(\,(197;84^{\times 4}, 42^{\times 6}) 
\#_{_{4}} u) \#_{_{3}} u) \#_{_{2}} u) \#_{_{1}}u,$$ with 
$u=42(2; 1^{\times 4})=(84; 42^{\times 4}).$ Applying the preceding process 
successively to $(197;84^{\times 4}, 42^{\times 6})$ we obtain again the 
vector $(1; 0^{\times 10})$ which is nef.

\section{The limits of the algorithm} \cntrs
\label{sect-Remarks}

In its present version, the algorithm described in Section~\ref{subsect-Alg} 
has the disadvantage that it does not tell which decomposition 
$v'_j=v_{j+1}\# u_{j+1}$ one should choose at each stage in order the whole 
process to end successfully. We would like to emphasize that this decision is  
sometimes crucial as the following example shows: \\
Let $v_0=(10; 3^{\times 11}) \in H_{11}$. if one tries to decompose $v_0$ as 
$v_0=(10;3^{\times 2}, 9) \# (9;3^{\times 9})$ the process will fail to 
give any information on $v_0$. The reason is that although $(9;3^{\times 9})$ 
is nef $(10;3^{\times 2}, 9)$ is not, and so we cannot apply the gluing 
theorem. However, the decomposition $v_0=(10;3^{\times 7}, 6) \# (6; 3^4)$ 
will eventually lead to a successful ending of the algorithm, thus proving 
that $v_0$ is nef.
It would be useful of course to find a rule for choosing the ``optimal'' 
decomposition at each stage.

Finally, let us mention one simple example for which it seems that the 
algorithm fails to give information always. Consider the vector 
$v_0=(19,6^{\times 10})$ which by Nagata's conjecture should be ample. 
However, it seems that the vector $v_0$ is {\em indecomposable} in the sense 
that it is impossible to find even nef vectors 
$v_1 \in H_{n_1}, u_1\in H_{k_1}$ with $n_1,k_1 < 10$, 
such that $v_0'=v_1\# v_2$ lies in the same orbit as $v_0$ under the Cremona 
action.
It would be interesting to find the precise conditions for an ample 
(resp. nef) vector $v$ to be indecomposable.

\section{Symplectic interpretations} \cntrs 
\label{sect-Symplectic} 

The purpose of this section is to explain the intuition which give rise 
to the gluing Theorem~\ref{thm-glue2}. Interestingly enough this comes from 
symplectic geometry.

Symplectic geometry is the branch of geometry dealing with the structure of 
symplectic manifolds which are by definition pairs, $(M,\Om)$, 
consisting of a smooth manifold $M$ and a non-degenerated closed differential 
2-form $\Om$. The reader is referred to~\cite{A-G} and~\cite{M-S} for the 
foundations.

Due to developments in this field of research in the last decade, many 
analogies has been discovered between symplectic and complex manifolds.
These become especially striking in dimension 4, 
where symplectic 4-manifolds play the role of complex surfaces.
In several cases it turned out that algebro-geometric considerations, remain 
true when properly translated into the symplectic category, and so gave rise  
to new theorems in the symplectic framework. This principle is reflected very 
well in the classification of rational and ruled symplectic manifolds of 
Lalonde and McDuff, in the symplectic packing theorems of 
McDuff and Polterovich, in Ruan's symplectization of the extremal rays theory 
etc.

In this paper we have, in some sense, reversed this direction of reasoning.
Our main theorem is in fact an algebro-geometric translation of a very simple 
symplectic fact arising from the theory of symplectic packing.
We refer the reader to~\cite{M-P} for an excellent exposition on the 
symplectic packing problem.

Recall from~\cite{M-P} that a {\em symplectic packing} of $(M,\Om)$ by $N$ 
balls of radii $\lam_1,\ldots, \lam_N$ is a symplectic embedding 
$$\phipack,$$ where $B(\lam_j)$ stands for the standard Euclidean closed ball 
of radius $\lam_j$ of the same dimension as $M$, endowed with its standard 
symplectic structure $\omstd=dx_1\wedge dy_1 + \ldots + dx_n \wedge dy_n$.

It was discovered by McDuff that every symplectic packing gives rise to a 
symplectic form $\wtldOm$ on the blow-up $\Th:\wtldM \rarw M$ of $M$ at the 
points $p_1=\vphi_1(0), \ldots p_N=\vphi_N(0)$. This form lies in the 
cohomology class 
$$[\wtldOm]=[\Th^* \Om]-\pi \sum_{j=1}^N \lam_j^2 e_j, \eqno(1)$$ 
where $e_j$ denotes the \PoD to the homology class of the exceptional 
divisor $E_j$ of the blow-up. 
This procedure is called {\em symplectic blowing-up}.

Conversely, given a symplectic form $\wtldOm$ on $\wtldM$ which is 
non-degenerated on the exceptional divisors $E_j$ and with cohomology class 
as in (1) above, one can perform {\em symplectic blowing-down} at the 
exceptional divisors and obtain a symplectic form $\Om$ on $M$ and a 
symplectic packing $\phipack$.

Consider the symplectic manifold $(\CC P^2, \sig)$ where $\sig$ is the 
Fubini-Studi \Khlr form normalized such that the area of a projective line 
is $\pi$. Its cohomology class is $\pi l$, where $l\in H^2(\CPTU, \ZZ)$ is 
the standard positive generator. 

Call a vector of positive numbers $(d;m_1,\ldots, m_k)$ {\em symplectic} if 
the cohomology class $$d \Th_V^*l - \sum_{j=1}^k m_j e_j$$ can be represented 
by a symplectic form $\wtld{\om}$ on some blow-up 
$\Th_V:V \rarw \CC P^2$ of $\CC P^2$ at some $k$ distinct points, 
in such a way that $\wtld{\om}$ is non-degenerated on the exceptional 
divisors.

Now, let $M$ be a complex surface and $\Th_p:\wtldM_p \rarw M$ its  
blow-up at the point $p\in M$ with exceptional divisor $E$. 
Denote by $e$ the \PoD to the homology class of $E$.

\begin{prop}
\label{prop-Symplectic}

Let $a \in H^2(M)$ and suppose that there exists a positive 
number $m$ such that the cohomology class $\Th_p^*a -me\in H^2(\wtldM_p)$ 
can be represented by a symplectic form whose restriction to $E$ is 
non-degenerated.\footnote{This means that $E$ is a symplectic submanifold 
with respect to this form.}
Then, for every symplectic vector $(m; \al_1, \ldots, \al_k)$, the 
cohomology class $\Th^*a-\sum_{j=1}^k \al_j e_j$ on the blow-up 
$\Th:\wtldM \rarw M$ at some $k$ points can be represented by a 
symplectic form.

\end{prop}

The proof is based on the following very simple observation. 
If $\Th_p^*a -me$ has a symplectic representative $\wtldOm$, then by 
symplectic blowing down one obtains a symplectic form $\Om$ on $M$ and an 
embedding $\vphi$ of a standard 4-dimensional ball of radius 
$\sqrt{\frac{m}{\pi}}$ into $(M, \Om)$. 
The same argument with slight modifications, applied to the vector 
$(m; \al_1, \ldots, \al_k)$, implies that the standard ball of radius 
$\sqrt{\frac{m}{\pi}}$ admits a symplectic packing, say $\phi$, by $k$ balls 
of radii $\sqrt{\frac{\al_1}{\pi}}, \ldots, \sqrt{\frac{\al_k}{\pi}}$.
Composing these two embeddings we conclude that $(M, \Om)$ admits a 
symplectic packing $\vphi \circ \phi$ of $k$ balls of radii 
$\sqrt{\frac{\al_1}{\pi}}, \ldots, \sqrt{\frac{\al_k}{\pi}}$. The 
proposition follows now from symplectic blowing-up.
For completeness, here are the precise arguments of the proof.

\pf 
Let $\wtldOm$ be a symplectic form on $\wtldM_p$ lying in the cohomology 
class $\Th_p^*a -me$ and suppose that the restriction of $\wtldOm$ to $E$ 
is non-degenerated. Applying symplectic blowing-down to $\wtldOm$ we obtain a 
symplectic form $\Om$ on $M$ lying in the cohomology class \ $a$ \ and a 
symplectic embedding $\vphi:B(\sqrt{\frac{m}{\pi}}) \rarw (M, \Om)$.

Let $\wtld{\om}$ be a symplectic form on the blow-up $\pi:V \rarw \CPTU$ of 
$\CPTU$ lying in the cohomology class $m\Th_V^*l-\sum_{j=1}^k \al_j e_j$ and 
whose restriction to the exceptional divisors $E_j$ is non-degenerated. 
Blowing-down symplectically we obtain a symplectic form $\om$ on $\CPTU$ 
lying in the cohomology class \ $ m l$ \ and a symplectic packing 
$\psi:\coprod_{j=1}^k B(\sqrt{\frac{\al_j}{\pi}}) \rarw (\CPTU, \om)$. 
Since any two cohomologous symplectic forms on $\CPTU$ are symplectomorphic 
we may assume that $\om=\frac{m}{\pi}\sig$. It can be proved by the methods 
of~\cite{M-P} that there exists a symplectic submanifold (with respect to 
$\om$) $L\subset M$, homologous to a projective line, which is disjoint from 
$\mbox{Image}\;\psi$. It is well known that 
$(\CPTU \stmin L, \frac{m}{\pi}\sig) \approx 
B(\sqrt{\frac{m}{\pi}})$. We thus obtain a symplectic packing 
$\phi:\coprod_{j=1}^k B(\sqrt{\frac{\al_j}{\pi}}) \rarw 
B(\sqrt{\frac{m}{\pi}})$. 

The composition $\vphi\circ\phi$ is a symplectic packing of $(M,\Om)$ by $k$ 
balls of radii $\sqrt{\frac{\al_1}{\pi}}, \ldots, \sqrt{\frac{\al_k}{\pi}}$.
Blowing-up symplectically with respect to this embedding yields a symplectic 
form on the blow-up $\Th: \wtldM \rarw M$ of $M$ at $k$ points, which 
lies in the cohomology class $\Th^* a - \sum_{j=1}^k \al_je_j$. \ \Qed 
\ \\

Let us try to translate Proposition~\ref{prop-Symplectic} to the language 
of algebraic geometry. Keeping in mind that in the symplectic category the 
role of \Khlr forms is played by symplectic forms and the role of 
complex submanifolds by symplectic submanifolds, the K\"{a}hlerian 
translation should read: {\em ``If the cohomology class $\Th_p^*a-me$ has a 
\Khlr representative than for every \Khlr vector $(m; \al_1, \ldots, \al_k)$ 
the cohomology class $\Th^*a - \sum_{j=1}^k \al_j e_j$ has a \Khlr 
representative too''.} \\ Here, we call a vector $(m; \al_1, \ldots, \al_k)$ 
\Khlr if the cohomology 
class $$m \Th_V^*l - \sum_{j=1}^k \al_j e_j$$ can be represented 
by a \Khlr form $\wtld{\om}$ on some simple rational surface 
$\Th_V:V \rarw \CC P^2$ obtained by blowing-up $\CC P^2$ at some $k$ distinct 
points.

Due to Lefschetz theorem on $(1,1)$ classes and Kodaira's embedding 
theorem it follows that on a complex manifold there is a bijection -- 
via Poincar\'{e} duality -- between the set of homology classes of ample 
$\QQ$-divisors and the set of rational cohomology classes which can be 
represented by \Khlr forms. Poincar\'{e} dualizing the 
``K\"{a}hlerian translation'' we are naturally led to the following: 
{\em \ \ ``Let $D$ be a divisor on $M$ such that $\Th_p^*D-mE$ is ample.
Then for every ample vector $(m; \al_1, \ldots, \al_k)$ the divisor 
$\Th^*D -\sum_{j=1}^k \al_j E_j$ is ample too''.} \ \ This is precisely the 
contents of Theorem~\ref{thm-glue2} for the case that $M$ is a simple 
rational surface. The technical machinery which made the whole translation 
rigorous is Shustin's curve gluing technique which we used in 
Section~\ref{sect-Gluing}.

We would like to emphasize that the same ``symplectic reasoning'' suggests 
that if we replace the surface $S$ in the statement of 
Theorem~\ref{thm-glue2} by any projective surface, the Theorem should remain 
true. Similar symplectic arguments suggest that an appropriate version 
of Theorem~\ref{thm-glue2} should hold also for higher dimensions than 2.
It would be of course interesting to know whether Theorem~\ref{thm-glue2} 
continues to hold for smooth algebraic surfaces over an arbitrary 
algebraically closed field. 
We leave these discussions to another opportunity.  

\subsection{Symplectic meaning of the remainders $\calR_N(\calL)$} \cntrsb 
\label{subsect-Symplectic_meaning}

In section~\ref{subsect-Applic} we have defined the homogeneous remainders 
$\calR_N(\calL)$ of an ample line bundle $\calL$ over a surface.
The definition naturally extends to $n$-dimensional smooth varieties $X$ 
in the following obvious way. Given $p_1, \ldots, p_N \in X$, set  
$$\calR(\calL, p_1,\ldots,p_N)=\frac{1}{\calL ^n}
\inf_{0\leq\eps\in \RR} \left\{ \calL_{\eps}^n \bigg| 
\calL_{\eps} = \pi^*\calL -\eps\sum_{j=1}^N E_j \quad \mbox{is nef} 
\right\},$$
where $\pi:\wtldX\rarw X$ is the blow-up of $X$ at the points
$p_1,\ldots,p_N$ with exceptional divisors $E_i=\pi^{-1}(p_i)$.
To get a more global invariant, define
$$\calR_N (\calL)= \inf \left\{ \calR(\calL, p_1,\ldots p_N) 
\mid p_1, \ldots p_N \in X \quad \mbox{are distinct points} \right\}.$$ 

Let us explain now the symplectic meaning of these constants. 
Let $(M, \Om)$ be a symplectic manifold. Following McDuff and Polterovich 
define the following quantity 
$$ v_N(M,\Om) = \sup_{\vphi, \,\lam} \frac{\mbox{Vol}(\mbox{Image\, }\vphi)}
{\mbox{Vol}(M,\Om)},$$
where $(\vphi, \lam)$ passes over all the possible symplectic packings 
$\vphi$ of $(M,\Om)$ with $N$ equal balls of varying radius $\lam$.
The volume of the manifolds is defined 
to be $\mbox{Vol}(M,\Om)=\int_M \frac{1}{n!} \Om^{\wedge n}$.

The constants $v_N(M, \Om)$ admit values between $0$ and $1$ and measure the 
maximal part of the volume of $(M, \Om)$ which can be filled by symplectic 
packing with $N$ equal balls. When $v_N=1$ we say that there exists
a {\em full packing} of $(M,\Om)$ by $N$ equal balls, while in the case
$v_N<1$ we say that there exists a {\em packing obstruction}.

In view of the preceding discussion it is easy to see that the homogeneous 
remainders $\calR_N(\calL)$ of an ample line bundle over a complex manifold 
$M$, play the algebro-geometric role of the quantity $1-v_N(M, \Om)$, 
where $\Om$ is a \Khlr form representing the first Chern class of 
$\calL, \; c_1(\calL)$.
In fact, it is not hard to prove that the following inequality holds
$$1-v_N(M,\Om) \leq \calR_N(\calL). \eqno(2)$$

Note that there are cases in which one always has equality in (2).
For example, it follows from the work of 
McDuff and Polterovich (see~\cite{M-P}) that this is the case for $\CPTU$ 
when $N < 9$. The point is that the symplectic cone and the 
\Khlr cone of del Pezzo surfaces coincide. 
Note that in (real) dimension 4, more is known about the constants 
$v_N$ than about $\calR_N$ (see~\cite{Bi}). It would be interesting to know 
whether there exist cases in which a strict inequality occurs in (2). \\

Let us conclude by pointing out another interesting approach to bounding 
Seshadri constants via symplectic packing, 
due to Lazarsfeld (see~\cite{Laz}). The idea is that 
given a \Khlr form $\Om$ on a complex manifold and a symplectic packing 
$\vphi$ which is also holomorphic, the symplectic blow-up of $\Om$ 
associated to $\vphi$ will be K\"{a}hler. This situation happens when the 
associated \Khlr metric on the image of $\vphi$ is flat. 
Applying this to the case of a principally polarized abelian variety, 
Lazarsfeld obtains non-trivial estimates on Seshadri constatnts of the 
corresponding ample divisor. 

\subsection*{Acknowledgments} I am extremely grateful to Prof. Eugenii 
Shustin for explaining me his approach to Viro method with singularities 
and for numerous enlightening discussions from which I learned a good deal. 
I would like to thank Prof. Leonid Polterovich for his interest in this work 
and for encouraging and interesting me to work on this project. 
Special thanks are due to Prof. Joseph Bernstein for the encouragement and 
for drawing my attention to interesting points which I was not aware of.
I have also benefited from discussions with Michael Thaddeus who gave me 
the reference to~\cite{Do-Or} and with Ilya Tyomkin who shared with me his 
insights on the problems discussed in this paper. 
I wish to thank these people.

\end{document}